\documentstyle[epsf]{article}

\newcommand{\bigzerou}{%
\smash{\lower1.7ex\hbox{\bg 0}}}
\newtheorem{theorem}{Theorem}
\newtheorem{prop}{Proposition}
\newtheorem{defi}{Definition}

\newtheorem{fact}{Fact}

\newtheorem{Rem}{Remark}

\newcommand{\ba}{\begin{eqnarray}}
\newcommand{\ea}{\end{eqnarray}}
\newcommand{\no}{\nonumber}

\newcommand{\mapright}[1]{%
\smash{\mathop{%
\hbox to 1.0cm{\rightarrowfill}}\limits^{#1}}}
\newcommand{\mapleft}[1]{%
\smash{\mathop{%
\hbox to 1.3cm{\leftarrowfill}}\limits^{#1}}}

\begin{document}
\title{
\begin{flushright}
  \begin{minipage}[b]{5em}
    \normalsize
    ${}$      \\
  \end{minipage}
\end{flushright}
{\bf  Direct Proof of Mirror Theorem of Projective Hypersurfaces up to degree $3$ Rational Curves }}
\author{Masao Jinzenji\\
\\
\it Division of Mathematics, Graduate School of Science \\
\it Hokkaido University \\
\it  Kita-ku, Sapporo, 060-0810, Japan\\
{\it e-mail address: jin@math.sci.hokudai.ac.jp}}
\maketitle
\begin{abstract}
In this paper, we directly derive generalized mirror transformation of projective hypersurfaces up to degree $3$ genus $0$ Gromov-Witten invariants 
by comparing Kontsevich localization formula with residue integral representation of the virtual structure constants.
We can easily generalize our method for rational curves of arbitrary 
degree except for combinatorial complexities.
\end{abstract}
\section{Introduction}
In this paper, we prove generalized mirror transformation of genus $0$ Gromov-Witten invariants of degree $k$ hypersurface in 
$CP^{N-1}$ (we denote it by $M_{N}^{k}$) up to degree $3$ rational curves. For this purpose, we introduce here the virtual structure 
constants of $M_{N}^{k}$, that were first defined in our ancient work with A. Collino \cite{cj}.
\begin{defi}
The virtual structure constants $\tilde{L}_{n}^{N,k,d}$ $(d\leq 3,\;\;\tilde{L}_{n}^{N,k,d}\neq 0\;\; \mbox{only if}\;\; 0\leq
n\leq N-1-(N-k)d )$ are rational numbers defined by the   
initial condition and the recursive formula:
\begin{eqnarray}\sum_{n=0}^{k-1}\tilde{L}_{n}^{N,k,1}w^{n}&=&
k\cdot\prod_{j=1}^{k-1}(jw+(k-j)), 
\label{ini1}
\end{eqnarray}
\begin{eqnarray}
\tilde{L}^{N,k,1}_{n}&=&\tilde{L}^{N+1,k,1}_{n}, 
\label{rec1}
\end{eqnarray}
\begin{eqnarray}
\frac{\tilde{L}^{N,k,2}_{n}}{2}&=&\frac{1}{2}\cdot\frac{\tilde{L}^{N+1,k,2}_{n-1}}{2}+
\frac{1}{2}\cdot\frac{\tilde{L}^{N+1,k,2}_{n}}{2}
+\frac{1}{2}\tilde{L}^{N+1,k,1}_{n}\cdot \tilde{L}^{N+1,k,1}_{n+(N-k)},
\label{rec2}
\end{eqnarray}
\begin{eqnarray}
\frac{\tilde{L}^{N,k,3}_{n}}{3}&=&\frac{2}{9}\cdot\frac{\tilde{L}^{N+1,k,3}_{n-2}}{3}+\frac{5}{9}
\cdot\frac{\tilde{L}^{N+1,k,3}_{n-1}}{3}
+\frac{2}{9}\cdot\frac{\tilde{L}^{N+1,k,3}_{n}}{3}\nonumber\\
&&+\frac{4}{9}\cdot\frac{\tilde{L}^{N+1,k,2}_{n-1}}{2}\cdot \tilde{L}^{N+1,k,1}_{n+2(N-k)}
+\frac{1}{3}\cdot\frac{\tilde{L}^{N+1,k,2}_{n}}{2}\cdot \tilde{L}^{N+1,k,1}_{n+2(N-k)}\nonumber\\
&&+\frac{2}{9}\cdot\frac{\tilde{L}^{N+1,k,2}_{n}}{2}\cdot \tilde{L}^{N+1,k,1}_{n+1+2(N-k)}\nonumber\\
&&+\frac{2}{9}\cdot\tilde{L}^{N+1,k,1}_{n-1}\cdot \frac{\tilde{L}^{N+1,k,2}_{n-1+(N-k)}}{2}
+\frac{1}{3}\cdot\tilde{L}^{N+1,k,1}_{n}\cdot \frac{\tilde{L}^{N+1,k,2}_{n-1+(N-k)}}{2}\nonumber\\
&&+\frac{4}{9}\cdot\tilde{L}^{N+1,k,1}_{n}\cdot \frac{\tilde{L}^{N+1,k,2}_{n+(N-k)}}{2}\nonumber\\
&&+\frac{1}{3}\cdot\tilde{L}^{N+1,k,1}_{n}\cdot \tilde{L}^{N+1,k,1}_{n+(N-k)}\cdot
\tilde{L}^{N+1,k,1}_{n+2(N-k)}.
\label{rec3}
\end{eqnarray}
\end{defi} 
In mirror computation of Gromov-Witten invariants of $M_{N}^{k}$, 
$\tilde{L}_n^{N,k,d}$ is used as the B-model analogue of 3-point Gromov-Witten invariant:
\begin{eqnarray}
&& \frac{1}{k}\langle{\cal O}_{h^{N-2-n}} {\cal O}_{h^{n-1+(N-k)d}}{\cal O}_{h} \rangle_{0,d}=
\frac{d}{k}\langle{\cal O}_{h^{N-2-n}} {\cal O}_{h^{n-1+(N-k)d}}\rangle_{0,d}=
\no\\
&&=\frac{d}{k}\int_{\overline{M}_{0,2}(CP^{N-1},d)}c_{top}\bigl(R^{0}(\pi_{*}ev_{3}^{*}({\cal O}(k)))\bigr)\wedge ev_{1}^{*}(h^{N-2-n})
\wedge ev_{2}^{*}(h^{n-1+(N-k)d}).
\label{corr}
\end{eqnarray}
In (\ref{corr}), $h$ is  hyperplane class of $CP^{N-1}$, $ \overline{M}_{0,n}(CP^{N-1},d) $ represents moduli space of degree $d$ 
stable maps from genus $0$ stable curve to $CP^{N-1}$ with $n$ marked points,  
$ev_{i}: \overline{M}_{0,n}(CP^{N-1},d) \rightarrow CP^{N-1}$ is evaluation map of  the $i$-th marked point and $\pi:\overline{M}_{0,3}(CP^{N-1},d) 
\rightarrow \overline{M}_{0,2}(CP^{N-1},d)$ is forgetful map.
Definition of the virtual structure constants for arbitrary degree $d\;(\geq1)$ can be seen in \cite{gm}. 
In our previous paper \cite{vs}, we conjectured a residue integral representation of $\tilde{L}_{n}^{N,k,d}$, 
which can be interpreted as a result of localization computation on the moduli space of 
Gauged Sigma Model. 
In the following, we prepare some notations to describe 
the formula we conjectured. First, we define rational functions in $u,v$ by,
\begin{eqnarray}
e(k,d;u,v)&:=& \prod_{m=0}^{kd}\bigl(\frac{mu+(kd-m)v}{d}\bigr),\no\\
t(N,d;u,v)&:=&\prod_{m=1}^{d-1}\bigl(\frac{mu+(d-m)v}{d}\bigr)^N.
\end{eqnarray}
Next, we introduce ordered partition of positive integer $d$: 
\begin{defi}
Let $OP_{d}$ be the set of ordered partitions of positive integer $d$:
\begin{equation}
OP_{d}=\{\sigma_{d}=(d_{1},d_{2},\cdots,d_{l(\sigma_{d})})\;\;|\;\;
\sum_{j=1}^{l(\sigma_{d})}d_{j}=d\;\;,\;\;d_{j}\in{\bf N}\}.
\label{part} 
\end{equation}
From now on, we denote a ordered partition $\sigma_{d}$ by 
$(d_{1},d_{2},\cdots,d_{l(\sigma_{d})})$. In (\ref{part}), we denote 
the length of the ordered partition $\sigma_{d}$ by $l(\sigma_{d})$.
\end{defi}
With these set up, the residue integral representation is given as follows: 
\begin{eqnarray}
\frac{\tilde{L}_{n}^{N,k,d}}{d}&=&\frac{1}{k}\sum_{\sigma_{d}\in OP_{d}}
\frac{1}{(2\pi\sqrt{-1})^{l(\sigma_{d})+1}\prod_{j=1}^{l(\sigma_{d})}d_{j}}\oint_{C_{0}} dx_{l(\sigma_d)}\cdots
\oint_{C_{0}} dx_0 (x_{0})^{N-2-n}(x_{l(\sigma_{d})})^{n-1+(N-k)d}\times\no\\
&&\times\prod_{j=0}^{l(\sigma_{d})}\frac{1}{(x_j)^N}
\prod_{j=1}^{l(\sigma_{d})-1}\frac{1}{kx_j\biggl( \frac{x_{j}-x_{j-1}}{d_{j}}+\frac{x_{j}-x_{j+1}}{d_{j+1}}\biggr)}
\prod_{j=1}^{l(\sigma_{d})}\frac{e(k,d_j;x_{j-1},x_{j})}{t(N,d_j;x_{j-1},x_j)}.
\label{int}
\end{eqnarray}
In (\ref{int}), $\frac{1}{2\pi\sqrt{-1}}\oint_{C_{0}}dx_{j}$ represents operation of taking residue at $x_{j}=0$. 
We have to mention here that the residue integral in (\ref{int}) severely depends on order of integration. 
To be more precise, we have to take the residues of 
$x_{j}$'s in descending (or ascending) order of subscript $j$. 
In the appendix of this paper, we prove,
\begin{theorem}
(\ref{int}) holds true if $d\leq3$.
\end{theorem}
We can indeed prove (\ref{int}) holds true for arbitrary $d$, but we write down proof for $d\leq3$ cases in this paper, mainly for 
convenience of space. Full proof will appear elsewhere. 
If $N\leq k$, the Gromov-Witten invariant $\frac{1}{k}\langle{\cal O}_{h^{N-2-n}} {\cal O}_{h^{n-1+(N-k)d}}{\cal O}_{h} \rangle_{0,d}$
and $\tilde{L}_n^{N,k,d}$ are different. In this case, we can write the former by weighted homogeneous polynomial of 
$\tilde{L}_m^{N,k,d'}\;(d'\leq d)$.
This formula is the generalized mirror transformation in our sense.   
Main result of this paper is a proof of this transformation up to $d=3$ case, that was conjectured in \cite{gt1}:
\begin{theorem}
\begin{eqnarray}
\frac{1}{k} \langle{\cal O}_{h^{N-2-n}}{\cal O}_{h^{n-1+N-k}}\rangle_{0,1} &=&\tilde{L}_{n}^{N,k,1}-\tilde{L}_{1+(k-N)}^{N,k,1},
\label{th1}
\end{eqnarray}
\begin{eqnarray}
\frac{1}{k}\langle{\cal O}_{h^{N-2-n}}{\cal O}_{h^{n-1+2(N-k)}}\rangle_{0,2}&=&\frac{1}{2}(\tilde{L}_{n}^{N,k,2}-\tilde{L}_{1+2(k-N)}^{N,k,2})
-\tilde{L}_{1+(k-N)}^{N,k,1}(\sum_{j=0}^{k-N}(\tilde{L}_{n-j}^{N,k,1} 
- \tilde{L}_{1+2(k-N)-j}^{N,k,1})),
\label{th2}
\end{eqnarray}
\begin{eqnarray}
\frac{1}{k}\langle{\cal O}_{h^{N-2-n}}{\cal O}_{h^{n-1+3(N-k)}}\rangle_{0,3} &=&\frac{1}{3}(\tilde{L}_{n}^{N,k,3}-\tilde{L}_{1+3(k-N)}^{N,k,3})
-\tilde{L}_{1+(k-N)}^{N,k,1}(\sum_{j=0}^{k-N}
(\tilde{L}_{n-j}^{N,k,2}-\tilde{L}_{1+3(k-N)-j}^{N,k,2})
+C_{1,1}^{N,k,3}(n))\nonumber\\
&&-\frac{1}{2}\tilde{L}_{1+2(k-N)}^{N,k,2}
(\sum_{j=0}^{2(k-N)}(\tilde{L}_{n-j}^{N,k,1}
-\tilde{L}_{1+3(k-N)-j}^{N,k,1}))\no\\
&&+\frac{3}{2}(\tilde{L}_{1+(k-N)}^{N,k,1})^{2}
(\sum_{j=0}^{2(k-N)}A_{j}(\tilde{L}_{n-j}^{N,k,1}
-\tilde{L}_{1+3(k-N)-j}^{N,k,1})),
\label{th3}
\end{eqnarray}
where
\begin{eqnarray}
A_{j}&:=&j+1,\;\;\mbox{if}\;\;\;(0\leq j\leq k-N),\;\;\;
A_{j}:=1+2(k-N)-j,\;\;\mbox{if}\;\;\; (k-N\leq j\leq 2(k-N)),\no\\
C_{1,1}^{N,k,3}(n)&=&\sum_{j=0}^{(k-N)-1}\bigl(\sum_{m=0}^{j}\tilde{L}_{n-m}^{N,k,1}
\tilde{L}_{n-2(k-N)+j-m}^{N,k,1}-\tilde{L}_{(k-N)+2+j}^{N,k,1}
(\sum_{m=0}^{2(k-N)}\tilde{L}_{n-m}^{N,k,1})\no\\
&&+\tilde{L}_{1+(k-N)}^{N,k,1}
(\sum_{m=j+1}^{2(k-N)-j-1}\tilde{L}_{n-m}^{N,k,1})\bigr)\no\\
&&-\sum_{j=0}^{(k-N)-1}\bigl(\sum_{m=0}^{j}\tilde{L}_{1+3(k-N)-m}^{N,k,1}
\tilde{L}_{1+(k-N)+j-m}^{N,k,1}-\tilde{L}_{(k-N)+2+j}^{N,k,1}
(\sum_{m=0}^{2(k-N)}\tilde{L}_{1+3(k-N)-m}^{N,k,1})\no\\
&&+\tilde{L}_{1+(k-N)}^{N,k,1}
(\sum_{m=j+1}^{2(k-N)-j-1}\tilde{L}_{1+3(k-N)-m}^{N,k,1})\bigr).
\label{fini}
\end{eqnarray}
\end{theorem}
Of course, the above formulas can be derived by using known methods presented in 
various papers: \cite{givc},\cite{iri},\cite{gene},\cite{blly}. 
In these works, generalized mirror transformation is derived as the effect of coordinate change of 
the B-model deformation parameters into the A-model ones.
We feel that this process  is a little bit too sophisticated to capture  
geometrical image of generalized mirror transformation: change of the moduli space of Gauged Linear Sigma Model 
into the one of stable maps. 
In this paper, we present an 
elementary and direct proof of Theorem 2 by using the result of Kontsevich \cite{kont} and Theorem 1. Our strategy is the following.
First, we write down explicit formula of $\langle{\cal O}_{h^{N-2-n}}{\cal O}_{h^{n-1+(N-k)d}}\rangle_{d}$ that follows 
from localization computation of Kontsevich. This formula includes combinatorially complicated summations with characters of 
torus action $\lambda_{j},\;(j=1,\cdots,N)$, but we can rewrite these summations into residue integrals of finite complex variables.  
This process is a generalization of the well-known computation on Bott residue theorem, that can be seen p.434-435 of \cite{gh}.
After this operation, we take non-equivariant limit $\lambda_{j}\rightarrow 0$. Resulting formula is very close to 
our residue integral representation of $\tilde{L}_{n}^{N,k,d}$ in (\ref{int}). With this formula, what we need for proof of Theorem 2 
is elementary combinatorial decomposition of rational functions in the integrands. 

This paper is organized as follows. In Section 2, we explain the process to reduce combinatorial summations in 
Kontsevich's localization formula to residue 
integrals in finite variables. Then we present residue integral representation of 2-point Gromov-Witten invariants. This representation can be 
directly compared with the r.h.s of (\ref{int}) after taking non-equivariant limit $\lambda_{j}\rightarrow 0$. In Section 3, we prove Theorem 2
by using decomposition of rational functions in the integrands. Section 3 gives concluding remarks. In Appendix, we 
prove Theorem 1, that plays important role in the proof of Theorem 2.    
\\
\\  
{\bf Acknowledgment} We would like to thank Dr. Brian Forbes for valuable discussions. 
We would also like to thank Miruko Jinzenji for kind encouragement. 
\section{Reduction of Localization Formula to Residue Integral}
We start from the Kontsevich's localization formulas for $2$-point genus 0 Gromov-Witten invariants of $M_{N}^{k}$. 
For compact presentation of these formulas, we introduce several notations: 
\begin{eqnarray}
w_{a}(u,v)&:=&\frac{u^{a}-v^{a}}{u-v}=\sum_{p+q=a-1,p,q\geq 0}u^p v^q,\no\\
w_{a}(u,v,w)&:=&\sum_{p+q+r=a-2,p,q,r\geq 0}u^p v^q w^r,\no\\
\end{eqnarray}
and, 
\begin{eqnarray}
E(k,d;i,j)&:=& \prod_{m=0}^{kd}\bigl(\frac{m\lambda_i+(kd-m)\lambda_j}{d}\bigr),\no\\
V(N;i)&:=&\prod_{j\neq i, 1\leq j\leq N}(\lambda_j-\lambda_i),\no\\
T(N,d;i,j)&:=&\prod_{k=1}^{N}\prod_{m=1}^{d-1}\bigl(\frac{m\lambda_i+(d-m)\lambda_j}{d}-\lambda_k\bigr).
\label{fun1}
\end{eqnarray}
In (\ref{fun1}), $\lambda_{j}\;(j=1,\cdots,N)$ are characters of torus action on $CP^{N-1}$:
\begin{equation}
(X_{1}:\cdots:X_{N}) \rightarrow (e^{\lambda_{1}t}X_{1}:\cdots:e^{\lambda_{N}t}X_{N}).  
\end{equation}
Here, we also introduce an elementary equality that will be used later in this paper:
\begin{equation}
w_{a}(x_1,x_2)+w_a(x_2,x_3)=(2x_2-x_1-x_3)w_a(x_1,x_2,x_3)+2w_a(x_1,x_3).
\label{rel1}
\end{equation}
With these set-up's, the localization formulas that represent $\langle {\cal O}_{h^a}{\cal O}_{h^b}\rangle_{0,d}\;(a=N-2-n,\;b=n-1+(N-k)d)$
 are described as follows:
\begin{fact}{\bf (Kontsevich)}
\begin{eqnarray}
\langle{\cal O}_{h^{a}}{\cal O}_{h^{b}}\rangle_{0,1}&=& -\frac{1}{2}\sum_{i\neq j}\frac{E(k,1;i,j)(\lambda_i-\lambda_j)^2}{V(N;i)V(N;j)}\cdot
w_a(\lambda_i,\lambda_j) w_b(\lambda_i,\lambda_j) ,  \label{fac1}
\end{eqnarray}
\begin{eqnarray}
\langle{\cal O}_{h^{a}}{\cal O}_{h^{b}}\rangle_{0,2}&=& -\frac{1}{4}\sum_{i\neq j}\frac{E(k,2;i,j)(\lambda_i-\lambda_j)^2}{T(N,2;i,j)V(N;i)V(N;j)}\cdot
w_a(\lambda_i,\lambda_j) w_b(\lambda_i,\lambda_j) +  \no\\
&&+\frac{1}{2}\sum_{i\neq j\neq l} \frac{E(k,1;i,j)E(k,1;j,l)}{V(N;i)V(N;j)V(N;l)k\lambda_j}\cdot
\frac{1}{\frac{1}{\lambda_j-\lambda_i}+\frac{1}{\lambda_j-\lambda_l}}\times\no\\
&&\times \bigl( w_a(\lambda_i,\lambda_j) + w_a(\lambda_j,\lambda_l) \bigr)
 \bigl( w_b(\lambda_i,\lambda_j) + w_b(\lambda_j,\lambda_l) \bigr),
\label{fac2}
\end{eqnarray}
\begin{eqnarray}
\langle{\cal O}_{h^{a}}{\cal O}_{h^{b}}\rangle_{0,3}&=&-\frac{1}{6}\sum_{i\neq j}\frac{E(k,3;i,j)(\lambda_i-\lambda_j)^2}{T(N,3;i,j)V(N;i)V(N;j)}\cdot
w_a(\lambda_i,\lambda_j) w_b(\lambda_i,\lambda_j)+  \no\\
&&+\frac{1}{2}\sum_{i\neq j\neq l}@
\frac{E(k,2;i,j)E(k,1;j,l)}{T(N,2;i,j)V(N;i)V(N;j)V(N;l)k\lambda_j}\cdot
\frac{1}{\frac{2}{\lambda_j-\lambda_i}+\frac{1}{\lambda_j-\lambda_l}} \times\no\\
&&\times \bigl(  2 w_a(\lambda_i,\lambda_j) + w_a(\lambda_j,\lambda_l)  \bigr)
 \bigl( 2 w_b(\lambda_i,\lambda_j) + w_b(\lambda_j,\lambda_l)  \bigr)-\no\\
&&-\frac{1}{2}\sum_{i\neq j\neq l\neq m} \frac{E(k,1;i,j)E(k,1;j,l)E(k,1;l,m)}{ V(N;i)V(N;j)V(N;l)V(N;m)k\lambda_j k\lambda_l }\times\no\\
&&\times\frac{1}{\frac{1}{\lambda_j-\lambda_i}+\frac{1}{\lambda_j-\lambda_l}}
\frac{1}{\frac{1}{\lambda_l-\lambda_j}+\frac{1}{\lambda_l-\lambda_m}}\cdot\frac{1}{(\lambda_j-\lambda_l)^2 } \times\no\\
&&\times \bigl(  w_a(\lambda_i,\lambda_j) + w_a(\lambda_j,\lambda_l) + w_a(\lambda_l,\lambda_m)   \bigr)\times\no\\
&&\times \bigl(  w_b(\lambda_i,\lambda_j) + w_b(\lambda_j,\lambda_l) + w_b(\lambda_l,\lambda_m)   \bigr)-\no\\
&&-\frac{1}{6}\sum_{i\neq j,i\neq l,i\neq m}\frac{E(k,1;i,j)E(k,1;i,l)E(k,1;i,m)}{V(N;i)V(N;j)V(N;l)V(N;m)(k\lambda_i)^2}\times \no\\
&&\times \bigl( w_a(\lambda_i,\lambda_j) + w_a(\lambda_i,\lambda_l) + w_a(\lambda_i,\lambda_m)   \bigr)\times\no\\
&&\times \bigl( w_b(\lambda_i,\lambda_j) + w_b(\lambda_i,\lambda_l) + w_b(\lambda_i,\lambda_m)   \bigr).
\label{fac3}
\end{eqnarray}
\end{fact}
\begin{Rem}
In (\ref{fac1}), (\ref{fac2}) and (\ref{fac3}), each summand corresponds to a tree graph that represents  degeneration type of stable maps \cite{kont}. The r.h.s.'s  of these 
equalities  are invariant under variation of characters of torus action, but in (\ref{fac2}) and (\ref{fac3}) each summand indeed varies under variation of them. 
\end{Rem}
Though some elementary simplification of complicated terms is operated, these formulas follow from the results in \cite{kont}. 
The above formulas includes many complicated summations, but we can rewrite these summations into residue integrals in finite complex variables as follows:  
\begin{prop}
\begin{eqnarray}
\langle{\cal O}_{h^{a}}{\cal O}_{h^{b}}\rangle_{1}&=& - \frac{1}{2} \frac{1}{(2\pi\sqrt{-1})^2} \oint_{C_0}dx_2\oint_{C_0}dx_1\frac{e(k,1;x_{1},x_{2})
(x_1-x_2)^2}{(x_1)^N(x_2)^N}\cdot w_a(x_1,x_2)w_b(x_1,x_2),  
\label{p1}
\end{eqnarray}
\begin{eqnarray}
\langle{\cal O}_{h^{a}}{\cal O}_{h^{b}}\rangle_{2}&=& -\frac{1}{4} \frac{1}{(2\pi\sqrt{-1})^2} \oint_{C_0}dx_2\oint_{C_0}dx_1
\frac{e(k,2;x_1,x_2)(x_1-x_2)^2}{(\frac{x_1+x_2}{2})^N(x_1)^N(x_2)^N}\cdot
w_a(x_1,x_2)w_b(x_1,x_2)+  \no\\
&&+ \frac{1}{2} \frac{1}{(2\pi\sqrt{-1})^3}\oint_{C_0}dx_3 \oint_{C_0}dx_2\oint_{C_0}dx_1
 \frac{e(k,1;x_1,x_2)e(k,1;x_2,x_3)}{(x_1)^N(x_2)^N(x_3)^N kx_2}\cdot
\frac{1}{\frac{1}{x_2-x_1}+\frac{1}{x_2-x_3}}\times\no\\
&&\times \bigl(w_a(x_1,x_2)+w_a(x_2,x_3)\bigr)
 \bigl(w_b(x_1,x_2)+w_b(x_2,x_3)\bigr),
\label{p2}
\end{eqnarray}
\begin{eqnarray}
\langle{\cal O}_{h^{a}}{\cal O}_{h^{b}}\rangle_{3}&=&-\frac{1}{6} \frac{1}{(2\pi\sqrt{-1})^2} \oint_{C_0}dx_2\oint_{C_0}dx_1
\frac{e(k,3;x_1,x_2)(x_1-x_2)^2}{(\frac{2x_1+x_2}{3})^N(\frac{x_1+2x_2}{3})^N (x_1)^N (x_2)^N}\cdot
w_a(x_1,x_2)w_b(x_1,x_2)+  \no\\
&&+\frac{1}{4} \frac{1}{(2\pi\sqrt{-1})^3}\oint_{C_0}dx_3 \oint_{C_0}dx_2\oint_{C_0}dx_1
\frac{e(k,2;x_1,x_2)e(k,1;x_2,x_3)}{(\frac{x_1+x_2}{2})^N (x_1)^N(x_2)^N(x_3)^N kx_2}\cdot
\frac{1}{\frac{2}{x_2-x_1}+\frac{1}{x_2-x_3}} \times\no\\
&&\times \bigl( 2w_a(x_1,x_2)+w_a(x_2,x_3)\bigr)
 \bigl( 2w_b(x_1,x_2)+w_b(x_2,x_3)\bigr)+\no\\
&&+\frac{1}{4} \frac{1}{(2\pi\sqrt{-1})^3}\oint_{C_0}dx_3 \oint_{C_0}dx_2\oint_{C_0}dx_1
\frac{e(k,1;x_1,x_2)e(k,2;x_2,x_3)}{(\frac{x_2+x_3}{2})^N (x_1)^N(x_2)^N(x_3)^N kx_2}\cdot
\frac{1}{\frac{1}{x_2-x_1}+\frac{2}{x_2-x_3}} \times\no\\
&&\times \bigl( w_a(x_1,x_2)+2w_a(x_2,x_3)\bigr) 
 \bigl( w_b(x_1,x_2)+2w_b(x_2,x_3)\bigr)-\no\\
&&-\frac{1}{2} \frac{1}{(2\pi\sqrt{-1})^4}\oint_{C_0}dx_4  \oint_{C_0}dx_3 \oint_{C_0}dx_2\oint_{C_0}dx_1
   \frac{e(k,1;x_1,x_2)e(k,1;x_2,x_3)e(k,1;x_3,x_4)}{(x_1)^N(x_2)^N(x_3)^N(x_4)^N kx_2 kx_3 }\times\no\\
&&\times\frac{1}{\frac{1}{x_2-x_1}+\frac{1}{x_2-x_3}}
\frac{1}{\frac{1}{x_3-x_2}+\frac{1}{x_3-x_4}}\cdot\frac{1}{(x_2-x_3)^2 } \times\no\\
&&\times \bigl(w_a(x_1,x_2)+w_a(x_2,x_3)+w_a(x_3,x_4) \bigr)\bigl(w_b(x_1,x_2)+w_b(x_2,x_3)+w_b(x_3,x_4) \bigr)-\no\\
&&-\frac{1}{6}\frac{1}{(2\pi\sqrt{-1})^4}\oint_{C_0}dx_4  \oint_{C_0}dx_3 \oint_{C_0}dx_2\oint_{C_0}dx_1
  \frac{e(k,1;x_1,x_2)e(k,1;x_1,x_3)e(k,1;x_1,x_4)}{(x_1)^N(x_2)^N(x_3)^N(x_4)^N(kx_1)^2}\times \no\\
&&\times \bigl(w_a(x_1,x_2)+w_a(x_1,x_3)+w_a(x_1,x_4) \bigr)\bigl(w_b(x_1,x_2)+w_b(x_1,x_3)+w_b(x_1,x_4) \bigr).
\label{p3}
\end{eqnarray}
\end{prop}
{\it proof)} For convenience of space, we write down proof of $d=1,2$ cases. Proof of $d=3$ case can be done in the same way 
as the proof below. 
We start from $d=1$ case. By elementary residue theorem, we can rewrite the r.h.s. of (\ref{fac1}) into the following residue integral:
\begin{eqnarray}
\langle{\cal O}_{h^{a}}{\cal O}_{h^{b}}\rangle_{1} &=& - \frac{1}{2} \frac{1}{(2\pi\sqrt{-1})^2} \oint_{C_{(0,5R)}}dx_2\oint_{C_{(0,R)}}dx_1\frac{e(k,1;x_{1},x_{2})
(x_1-x_2)^2}{\prod_{j=1}^{N}\bigl((x_1-\lambda_j)(x_2-\lambda_j)\bigr)}\cdot w_a(x_1,x_2)w_b(x_1,x_2).\no\\ 
\label{line}
\end{eqnarray}
In (\ref{line}), $C(0,a)$ denotes the circle with center at $0$ and with radius $a$, and $R$ is a sufficiently large positive real 
number greater than $\mbox{max.}\{3|\lambda_{j}|\;|\; j=1,\cdots,N\}$.
In $d=2$ case, we also have similar equalities:
\begin{eqnarray}
&& -\frac{1}{4} \frac{1}{(2\pi\sqrt{-1})^2} \oint_{C_{(0,5R)}}dx_2\oint_{C_{(0,R)}}dx_1
\frac{e(k,2;x_1,x_2)(x_1-x_2)^2}{\prod_{j=1}^N\bigl((\frac{x_1+x_2}{2}-\lambda_j)(x_1-\lambda_j)(x_2-\lambda_j)\bigr)}\cdot
w_a(x_1,x_2)w_b(x_1,x_2) =\no\\
&&=-\frac{1}{4} \sum_{i\neq j}
\frac{E(k,2;i,j)(\lambda_i-\lambda_j)^2}{T(N,2;i,j)V(N;i)V(N;j)}\cdot
w_a(\lambda_i,\lambda_j)w_b(\lambda_i,\lambda_j) -\no\\
&&-2\sum_{i\neq j}\frac{ e(k,2;\lambda_i,2\lambda_j-\lambda_i) (\lambda_i-\lambda_j)^2}{ \bigl(\prod_{l=1}^{N}(2\lambda_j-\lambda_i-\lambda_l)\bigr)
V(N;i)V(N;j)}\cdot
w_a(\lambda_i,2\lambda_j-\lambda_i)w_b(\lambda_i,2\lambda_j-\lambda_i),
\label{conic1} 
\end{eqnarray}
and,
\begin{eqnarray}
&& \frac{1}{2} \frac{1}{(2\pi\sqrt{-1})^3}\oint_{C_{(0,25R)}}dx_3 \oint_{C_{(0,5R)}}dx_2\oint_{C_{(0,R)}}dx_1
 \frac{e(k,1;x_1,x_2)e(k,1;x_2,x_3)}{kx_2\prod_{j=1}^N\bigl((x_1-\lambda_j)(x_2-\lambda_j)(x_3-\lambda_j) \bigr)}\times\no\\
&&\times\frac{(x_2-x_1)(x_2-x_3)}{2x_2-x_1-x_3}\bigl(w_a(x_1,x_2)+w_a(x_2,x_3)\bigr)
 \bigl(w_b(x_1,x_2)+w_b(x_2,x_3)\bigr) =\no\\
&&= \frac{1}{2} \sum_{i\neq j,j\neq l}
 \frac{E(k,1;i,j)E(k,1;j,l)}{k\lambda_j V(N;i)V(N;j)V(N;l)}\frac{(\lambda_j-\lambda_i)(\lambda_j-\lambda_l)}{2\lambda_j-\lambda_i-\lambda_l}
\bigl(w_a(\lambda_i,\lambda_j)+w_a(\lambda_j,\lambda_l)\bigr)
 \bigl(w_b(\lambda_i,\lambda_j)+w_b(\lambda_j,\lambda_l)\bigr) + \no\\
&&+\frac{1}{2}\sum_{i\neq j}\frac{e(k;1;\lambda_i,\lambda_j)e(k,1;\lambda_j,2\lambda_j-\lambda_i)(\lambda_i-\lambda_j)^2}
{k\lambda_j V(N;i)V(N;j)V(N;l)\bigl(\prod_{l=1}^{N}(2\lambda_j-\lambda_i-\lambda_l)\bigr)}
\bigl( w_a(\lambda_i,\lambda_j)+w_a(\lambda_j,2\lambda_j-\lambda_i) \bigr)\times\no\\
 &&\times\bigl(w_b(\lambda_i,\lambda_j)+w_b(\lambda_j,2\lambda_j-\lambda_i)\bigr).
\label{conic2} 
\end{eqnarray}
On the other hand, we can easily see the following relations: 
\begin{eqnarray}
&&e(k,2;\lambda_i,2\lambda_j-\lambda_i) =\frac{e(k;1;\lambda_i,\lambda_j)e(k,1;\lambda_j,2\lambda_j-\lambda_i)}{k\lambda_j},\no\\
&& w_a(\lambda_i,\lambda_j)+w_a(\lambda_j,2\lambda_j-\lambda_i) = 2w_a(\lambda_i,2\lambda_j-\lambda_i).
\label{rel2}
\end{eqnarray}
The second equality follows from (\ref{rel1}).
With these relations, the second terms of the r.h.s.'s of (\ref{conic1}) and (\ref{conic2}) cancel.  
Then we obtain,
\begin{eqnarray}
&&\langle{\cal O}_{h^{a}}{\cal O}_{h^{b}}\rangle_{2}=\no\\
&&=-\frac{1}{4} \frac{1}{(2\pi\sqrt{-1})^2} \oint_{C_{(0,5R)}}dx_2\oint_{C_{(0,R)}}dx_1
\frac{e(k,2;x_1,x_2)(x_1-x_2)^2}{\prod_{j=1}^N\bigl((\frac{x_1+x_2}{2}-\lambda_j)(x_1-\lambda_j)(x_2-\lambda_j)\bigr)}\cdot
w_a(x_1,x_2)w_b(x_1,x_2) +\no\\
&&+\frac{1}{2} \frac{1}{(2\pi\sqrt{-1})^3}\oint_{C_{(0,25R)}}dx_3 \oint_{C_{(0,5R)}}dx_2\oint_{C_{(0,R)}}dx_1
 \frac{e(k,1;x_1,x_2)e(k,1;x_2,x_3)}{kx_2\prod_{j=1}^N\bigl((x_1-\lambda_j)(x_2-\lambda_j)(x_3-\lambda_j) \bigr)}\times\no\\
&&\times\frac{(x_2-x_1)(x_2-x_3)}{2x_2-x_1-x_3}\bigl(w_a(x_1,x_2)+w_a(x_2,x_3)\bigr)
 \bigl(w_b(x_1,x_2)+w_b(x_2,x_3)\bigr).
\label{conic}
\end{eqnarray}
If we look at the r.h.s.'s of (\ref{line}) and (\ref{conic}), we can easily see by coordinate change $x_{j}=\frac{1}{z_{j}}$, that 
each summand is invariant under variation of $\lambda_{j}$'s. Therefore, we can take non-equivariant limit $\lambda_j\rightarrow 0$ (we also 
take $R\rightarrow 0$ limit). This operation leads us to the equalities of the proposition.  
\begin{Rem}
In Remark 1, we noticed that each summand in (\ref{fac2}) is not invariant under variation of 
characters. But as can be seen in (\ref{conic1}) and (\ref{conic2}), we can make it invariant by adding 
a suitable rational function of characters. These additional rational functions chancel out after 
adding up summands that correspond to tree graphs. The same mechanism also works in $d=3$ case.    
\end{Rem}
\section{Proof of Theorem 2}
Before we go into proof of Theorem 2, we note an equality: 
\begin{equation}
\tilde{L}_n^{N,k,d}=\tilde{L}_{N-1-(N-k)d-n}^{N,k,d},
\label{reflection}
\end{equation}
, which naturally follows from Theorem 1. 

Let us start from $d=1$ case. In this case, we apply a trivial equality:
\begin{eqnarray}
(x_1-x_2)^2w_a(x_1,x_2)w_b(x_1,x_2)=(x_1^a-x_2^a)(x_2^b-x_1^b)=x_1^a x_2^b+x_1^b x_2^a-x_1^{a+b}-x_2^{a+b},
\label{dec1}
\end{eqnarray}
to the r.h.s of (\ref{p1}). Then the theorem follows from Theorem 1 and (\ref{reflection}).
In $d=2$ case, we apply (\ref{dec1}) to the first summand of the r.h.s. of (\ref{p2}). To the 
second summand of it, we apply the following decomposition of the rational function in the integrand: 
\begin{eqnarray}
&&\frac{(w_a(x_1,x_2)+w_a(x_2,x_3))(w_b(x_1,x_2)+w_b(x_2,x_3))}{\frac{1}{x_2-x_1}+\frac{1}{x_2-x_3}}=\no\\
&&=\frac{(x_2-x_1)(x_2-x_3)(w_a(x_1,x_2)+w_a(x_2,x_3))(w_b(x_1,x_2)+w_b(x_2,x_3))}{2x_{2}-x_{1}-x_{3}}=\no\\
&&=\frac{((x_1)^a-(x_3)^a)((x_3)^b-(x_1)^b)}{2x_2-x_1-x_3}+(x_2-x_1)w_a(x_1,x_2)w_b(x_1,x_2)+
(x_2-x_3)w_a(x_2,x_3)w_b(x_2,x_3).\no\\
\label{dec2}
\end{eqnarray}
Then the first summand in the r.h.s. of  (\ref{p2}) and the first term in the decomposition (\ref{dec2}) add up to, 
\begin{equation}
\frac{1}{2}k(\tilde{L}_{n}^{N,k,d}-\tilde{L}_{1+2(k-N)}),
\end{equation}  
by using (\ref{reflection}) and Theorem 1. The second and the third terms in the decomposition (\ref{dec2}) result in,
\begin{equation}
k\sum_{j=0}^{b-1}\tilde{L}^{N,k,1}_{1+k-N}( \tilde{L}_{1+j+k-N}^{N,k,1}-\tilde{L}_{1+a+j+k-N}^{N,k,1}),
\label{gm2}
\end{equation}
by using $(x_i-x_j)w_a(x_i,x_j)=x_i^a-x_j^a$, (\ref{reflection}) and Theorem 1. But if we look back at $a=N-2-n,\; b=n-1-2(k-N)$ and (\ref{reflection}), 
(\ref{gm2}) turns out to be, 
\begin{equation}
-k\sum_{j=0}^{k-N}\tilde{L}^{N,k,1}_{1+k-N}( \tilde{L}_{n-j}^{N,k,1} - \tilde{L}_{1+2(k-N)-j}^{N,k,1} ).
\label{f2}
\end{equation}
This completes the proof of $d=2$ case.

Now, we turn into $d=3$ case. In the same way as the $d=1,2$ cases, we apply (\ref{dec1}) to the first summand of the r.h.s. of 
(\ref{p3}). To the second and the third summands of it, we apply the decompositions:
\begin{eqnarray}
&&\frac{2(2w_a(x_1,x_2)+w_a(x_2,x_3))(2w_b(x_1,x_2)+w_b(x_2,x_3))}{\frac{2}{x_2-x_1}+\frac{1}{x_2-x_3}}=\no\\
&&=\frac{(x_2-x_1)(x_2-x_3)(2w_a(x_1,x_2)+w_a(x_2,x_3))(2w_b(x_1,x_2)+w_b(x_2,x_3))}{\frac{x_{2}-x_{1}}{2}+x_2-x_{3}}=\no\\
&&=2\frac{((x_1)^a-(x_3)^a)((x_3)^b-(x_1)^b)}{\frac{x_2-x_1}{2}+x_2-x_3}+4(x_2-x_1)w_a(x_1,x_2)w_b(x_1,x_2)+
2(x_2-x_3)w_a(x_2,x_3)w_b(x_2,x_3).\no\\
\label{dec3}
\end{eqnarray}
and, 
\begin{eqnarray}
&&\frac{2(w_a(x_1,x_2)+2w_a(x_2,x_3))(w_b(x_1,x_2)+2w_b(x_2,x_3))}{\frac{1}{x_2-x_1}+\frac{2}{x_2-x_3}}=\no\\
&&=\frac{(x_2-x_1)(x_2-x_3)(w_a(x_1,x_2)+2w_a(x_2,x_3))(w_b(x_1,x_2)+2w_b(x_2,x_3))}{x_{2}-x_{1}+\frac{x_2-x_{3}}{2}}=\no\\
&&=2\frac{((x_1)^a-(x_3)^a)((x_3)^b-(x_1)^b)}{x_2-x_1+\frac{x_2-x_3}{2}}+2(x_2-x_1)w_a(x_1,x_2)w_b(x_1,x_2)+
4(x_2-x_3)w_a(x_2,x_3)w_b(x_2,x_3).\no\\
\label{dec4}
\end{eqnarray}
Lastly, we apply the following decomposition to the fourth summand.
\begin{eqnarray}
&&-\frac{(w_a(x_1,x_2)+w_a(x_2,x_3)+w_a(x_3,x_4))(w_b(x_1,x_2)+w_b(x_2,x_3)+w_b(x_3,x_4))}
{(\frac{1}{x_2-x_1}+\frac{1}{x_2-x_3})(\frac{1}{x_3-x_2}+\frac{1}{x_3-x_4})(x_2-x_3)^2}\no\\
&&=\frac{(x_2-x_1)(x_3-x_4)(w_a(x_1,x_2)+w_a(x_2,x_3)+w_a(x_3,x_4))(w_b(x_1,x_2)+w_b(x_2,x_3)+w_b(x_3,x_4))}{r_1r_2}\no\\
&&= \frac{((x_1)^a-(x_4)^a)((x_4)^b-(x_1)^b)}{r_1r_2}+\no\\
&& +\frac{2(x_3-x_1)w_a(x_1,x_3)w_b(x_1,x_3)+(x_3-x_4)w_a(x_3,x_4)w_b(x_3,x_4)}{r_1}+\no\\
&& +\frac{2(x_2-x_4)w_a(x_2,x_4)w_b(x_2,x_4)+(x_2-x_1)w_a(x_1,x_2)w_b(x_1,x_2)}{r_2}+\no\\
&&+\frac{1}{2}\bigl( w_a(x_1,x_2)w_b(x_1,x_2)+w_a(x_3,x_4)w_b(x_3,x_4)+w_a(x_1,x_2)w_b(x_2,x_3)+w_a(x_2,x_3)w_b(x_1,x_2)+\no\\
&&+w_a(x_2,x_3)w_b(x_3,x_4)+w_a(x_3,x_4)w_b(x_2,x_3)\bigr)+\no\\
&&+(x_{3}-x_{1})\bigl(w_a(x_1,x_3)w_b(x_1,x_2,x_3)+w_a(x_1,x_2,x_3)w_b(x_1,x_3)+\frac{1}{2}r_1 w_a(x_1,x_2,x_3)w_b(x_1,x_2,x_3)\bigr)+\no\\
&&+(x_{2}-x_{4})\bigl(w_a(x_2,x_4)w_b(x_2,x_3,x_4)+w_a(x_2,x_3,x_4)w_b(x_2,x_4)+\frac{1}{2}r_2 w_a(x_2,x_3,x_4)w_b(x_2,x_3,x_4)\bigr),
\label{dec5}
\end{eqnarray}
where $r_1=2x_2-x_1-x_3,\;r_2=2x_3-x_2-x_4$.
The first summand of the r.h.s. of (\ref{p3}), the first terms of (\ref{dec3}) and (\ref{dec4}), and the term with denominator $r_1 r_2$ in (\ref{dec5})
add up to $\frac{1}{3}k(\tilde{L}_n^{N,k,3}-\tilde{L}_{1+3(k-N)}^{N,k,3})$ by (\ref{reflection}) and Theorem 1.
The second term of (\ref{dec3}), the third term of (\ref{dec4}) and
\begin{equation}
\frac{2(x_3-x_1)w_a(x_1,x_3)w_b(x_1,x_3)}{r_1} ,\;\;\frac{2(x_2-x_4)w_a(x_2,x_4)w_b(x_2,x_4)}{r_2}, 
\end{equation}
in (\ref{dec5}) add up to,
 \begin{equation}
-k\sum_{j=0}^{k-N}\tilde{L}^{N,k,1}_{1+k-N}(\tilde{L}_{n-j}^{N,k,2} -  \tilde{L}_{1+3(k-N)-j}^{N,k,2} ).
\end{equation}
by Theorem 1 and the same reorganization used to derive (\ref{f2}).
The third term of (\ref{dec3}), the second term of (\ref{dec4}) and
\begin{equation}
\frac{(x_3-x_4)w_a(x_3,x_4)w_b(x_3,x_4)}{r_1} ,\;\;\frac{(x_2-x_1)w_a(x_1,x_2)w_b(x_1,x_2)}{r_2}, 
\end{equation}
in (\ref{dec5}) add up to,
 \begin{equation}
-k\cdot\frac{1}{2}\sum_{j=0}^{2(k-N)}\tilde{L}^{N,k,2}_{1+2(k-N)}( \tilde{L}_{n-j}^{N,k,1} - \tilde{L}_{1+3(k-N)-j}^{N,k,1} ).
\end{equation}
in the same way as the previous argument.
Remaining terms of the decomposition (\ref{dec5}) and the fifth summand of the r.h.s. of (\ref{p3}) are reorganized as follows:
\begin{eqnarray}
&& k\tilde{L}^{N,k,1}_{1+k-N}\biggl(\frac{1}{2} \sum_{i=0}^{a-1}\sum_{j=0}^{b-1} \tilde{L}^{N,k,1}_{1+i+j+k-N} \tilde{L}^{N,k,1}_{2+2(k-N)} +\no\\
&&+ \frac{1}{2} \sum_{i=0}^{a-1}\sum_{j=0}^{b-1} \tilde{L}^{N,k,1}_{1+i+k-N}\tilde{L}^{N,k,1}_{2+j+2(k-N)}
+ \frac{1}{2} \sum_{i=0}^{a-1}\sum_{j=0}^{b-1} \tilde{L}^{N,k,1}_{1+j+k-N}\tilde{L}^{N,k,1}_{2+i+2(k-N)}+ \no\\
&& + \frac{1}{2} \sum_{i=0}^{a-1}\sum_{j=0}^{b-2}\bigl( \tilde{L}^{N,k,1}_{1+i+k-N}\tilde{L}^{N,k,1}_{3+j+2(k-N)}-
 \tilde{L}^{N,k,1}_{2+i+j+k-N}\tilde{L}^{N,k,1}_{2+2(k-N)}\bigr)+ \no\\
&& + \frac{1}{2} \sum_{i=0}^{a-1}\sum_{j=0}^{b-2}\bigl( \tilde{L}^{N,k,1}_{1+k-N}\tilde{L}^{N,k,1}_{3+i+j+2(k-N)}-
 \tilde{L}^{N,k,1}_{2+j+k-N}\tilde{L}^{N,k,1}_{2+i+2(k-N)}\bigr)+\no\\
&& + \sum_{i=0}^{a-2}\sum_{j=0}^{i}\bigl( \tilde{L}^{N,k,1}_{1+j+k-N}\tilde{L}^{N,k,1}_{a+b-i+2(k-N)}-
 \tilde{L}^{N,k,1}_{1+b+j+k-N}\tilde{L}^{N,k,1}_{a-i+2(k-N)}\bigr)-\no\\
&& - \frac{1}{2} \sum_{i=0}^{a-1}\sum_{j=0}^{b-1} \tilde{L}^{N,k,1}_{1+k-N}\tilde{L}^{N,k,1}_{1+i+j+k-N}
-\sum_{i=0}^{a-1}\sum_{j=0}^{b-1} \tilde{L}^{N,k,1}_{1+i+k-N}\tilde{L}^{N,k,1}_{1+j+k-N}\biggr)=\no\\
&&= k\tilde{L}^{N,k,1}_{1+k-N}\biggl(-\sum_{i=0}^{a-1}\sum_{j=0}^{i} \bigl( \tilde{L}^{N,k,1}_{n+j-2(k-N)}\tilde{L}^{N,k,1}_{n+i+1-(k-N)} 
-\tilde{L}^{N,k,1}_{1+j+(k-N)}\tilde{L}^{N,k,1}_{2+i+2(k-N)}\bigr) +\no\\
&& +\sum_{i=0}^{a-1}\sum_{j=0}^{k-N} \tilde{L}^{N,k,1}_{n+1+i}\bigl(\tilde{L}^{N,k,1}_{n+j-2(k-N)}-
 \tilde{L}^{N,k,1}_{1+j+(k-N)}\bigr)+ \no\\
&& + \sum_{j=0}^{a-1}\bigl( \tilde{L}^{N,k,1}_{1+2(k-N)}\tilde{L}^{N,k,1}_{n+j-2(k-N)}-
 \tilde{L}^{N,k,1}_{n-(k-N)}\tilde{L}^{N,k,1}_{1+j+(k-N)}\bigr)+\no\\
&& + \frac{1}{2} \sum_{i=0}^{a-1}\sum_{j=0}^{k-N} \tilde{L}^{N,k,1}_{1+k-N}\bigl(\tilde{L}^{N,k,1}_{2+i+j+2(k-N)}
-\tilde{L}^{N,k,1}_{1+i+j+(k-N)}\bigr)\biggr)=\no\\
&&=k\tilde{L}^{N,k,1}_{1+k-N}\biggl(-\sum_{i=0}^{(k-N)-1}\sum_{j=0}^{i} \bigl( \tilde{L}^{N,k,1}_{n+j-2(k-N)}\tilde{L}^{N,k,1}_{n+i+1-(k-N)} 
-\tilde{L}^{N,k,1}_{1+j+(k-N)}\tilde{L}^{N,k,1}_{2+i+2(k-N)}\bigr) +\no\\
&& +\sum_{i=0}^{2(k-N)}\sum_{j=0}^{k-N} \tilde{L}^{N,k,1}_{1+j+(k-N)}\bigl(\tilde{L}^{N,k,1}_{n-i}-
 \tilde{L}^{N,k,1}_{1+i+(k-N)}\bigr)+ \no\\
&& + \frac{1}{2}\sum_{i=0}^{k-N}\sum_{j=0}^{k-N} \tilde{L}^{N,k,1}_{1+k-N}\bigl(\tilde{L}^{N,k,1}_{n-i-j}
-\tilde{L}^{N,k,1}_{1+i+j+(k-N)}\bigr)\biggr)=\no\\
&&=- k\tilde{L}^{N,k,1}_{1+k-N} C_{1,1}^{N,k,3}(n)+
k\cdot\frac{3}{2}(\tilde{L}^{N,k,1}_{1+k-N})^{2}(\sum_{j=0}^{2(k-N)}A_{j}(\tilde{L}_{n-j}^{N,k,1}
-\tilde{L}_{1+3(k-N)-j}^{N,k,1})). 
\end{eqnarray}
In this derivation, we only use the conditions:
\begin{equation}
a=N-2-n,\;\;b=n-1-3(k-N),\;\;\tilde{L}^{N,k,1}_{n}=\tilde{L}^{N,k,1}_{k-1-n},
\end{equation}
but need careful treatment of summations. Anyway, the final formula completes the proof of Theorem 2.   
\section{Conclusion}
Our motivation of this paper is to understand explicitly difference between the moduli space of Gauged 
Linear Sigma Model and the moduli space of stable maps. Unfortunately, this paper's treatment is quite 
computational. It gives us an elementary and explicit proof of mirror theorem for 
rational curves of lower degrees but lacks geometrical vision. We can indeed extend this paper's method 
to rational curves of higher degrees because we don't use geometrical simplicity of moduli space of 
rational curves of lower degrees. As can be seen in \cite{gene} and \cite{vs}, generalized mirror 
transformation of $d=4,5$ rational curves has quite complicated structure. Therefore, we need 
combinatorial sophistication of our method. For this purpose, we had better search for geometrical 
meaning of decomposition of rational functions, such as (\ref{dec5}). 

One of the main features of this paper may be translation of combinatorial summations in Kontsevich's 
localization formula into residue integrals, that enabled us to compare directly the Gromov-Witten 
invariants with the virtual structure constants. This translation can be applied to various examples. 
At least, we can use it to prove mirror theorem of ${\cal O}(1)\oplus {\cal O}(-3)\rightarrow {\bf P^{1}}$
\cite{fj1}. We also think that we can apply it to prove mirror theorem at higher genus. Anyway, we 
have to pursue combinatorial sophistication of our method.   
\newpage
\appendix
\section{Appendix: Proof of Theorem 1}
We prove Theorem 1 by showing that the r.h.s. of (\ref{int}) satisfies the initial condition and the 
recursive formulas (\ref{ini1}), (\ref{rec1}), (\ref{rec2}) and (\ref{rec3}).
For this purpose, we note here a relation between rational functions that appear in the residue integrals:
\begin{eqnarray}
&&\prod_{j=0}^{l(\sigma_{d})}\frac{1}{(x_j)^N}
\prod_{j=1}^{l(\sigma_{d})-1}\frac{1}{kx_j\biggl( \frac{x_{j}-x_{j-1}}{d_{j}}+\frac{x_{j}-x_{j+1}}{d_{j+1}}\biggr)}
\prod_{j=1}^{l(\sigma_{d})}\frac{e(k,d_j;x_{j-1},x_{j})}{t(N,d_j;x_{j-1},x_j)}=\no\\
&&=\prod_{j=0}^{l(\sigma_{d})}\frac{1}{(x_j)^{N+1}}
\prod_{j=1}^{l(\sigma_{d})-1}\frac{1}{kx_j\biggl( \frac{x_{j}-x_{j-1}}{d_{j}}+\frac{x_{j}-x_{j+1}}{d_{j+1}}\biggr)}
\prod_{j=1}^{l(\sigma_{d})}\frac{e(k,d_j;x_{j-1},x_{j})}{t(N+1,d_j;x_{j-1},x_j)}\times\no\\
&&\times (x_0 x_1\cdots x_{l(\sigma_{d})})\prod_{j=1}^{l(\sigma_{d})}\prod_{i=1}^{d_j-1}(\frac{ix_{j-1}+(d_j-i)x_{j}}{d_j}).
\label{reduction}
\end{eqnarray}
In $d=1$ case, the r.h.s of (\ref{int}) becomes,
\begin{eqnarray}
&& \frac{1}{k}\cdot\frac{1}{(2\pi\sqrt{-1})^2}\oint_{C_{0}}dx_1\oint_{C_{0}}dx_0\; x_0^{N-2-n}x_{1}^{n-1+N-k}\frac{e(k,1;x_0,x_1)}{x_0^{N}x_1^{N}} =\no\\
&&=k\frac{1}{(2\pi\sqrt{-1})^2}\oint_{C_{0}}dx_1\oint_{C_{0}}dx_0 \frac{\prod_{j=1}^{k-1}(jx_0+(k-j)x_1)}{x_0^{n+1}x_1^{k-n}}. 
\end{eqnarray}
Hence (\ref{ini1}) and (\ref{rec1}) automatically hold true.

(\ref{reduction}) tells us that the recursive formulas in $d=2,3$ cases follow from adequate decomposition of, 
$$ (x_0 x_1\cdots x_{l(\sigma_{d})})\prod_{j=1}^{l(\sigma_{d})}\prod_{i=1}^{d_j-1}(\frac{ix_{j-1}+(d_j-i)x_{j}}{d_j}).$$
Explicitly, decompositions are given as follows: \\
{\bf d=2 case}
\begin{eqnarray}
&& \sigma_{2}=(2):\;\;\;\; x_0x_1\frac{x_0+x_1}{2},\no\\
&& \sigma_{2}=(1,1):\;\;\; x_0x_1x_2=x_0x_2\frac{x_0+x_2}{2}+\frac{1}{2}rx_0x_2,\;\;(r=2x_1-x_0-x_2).
\end{eqnarray}
{\bf d=3 case}
\begin{eqnarray}
&&\sigma_{3}=(3):\;\;\;\; x_0x_1\frac{2x_0+x_1}{3}\frac{x_0+2x_1}{3}=x_0x_1(\frac{2}{9}x_0^2+\frac{5}{9}x_0x_1+\frac{2}{9}x_1^2),\no\\
&&\sigma_{3}=(2,1):\;\;\;\; x_0x_1x_2\frac{x_0+x_1}{2}=x_0x_2(\frac{2}{9}x_0^2+\frac{5}{9}x_0x_2+\frac{2}{9}x_2^2+
r_1(\frac{4}{9}x_0+\frac{1}{3}x_1+\frac{2}{9}x_2)),\no\\
&&(r_1=\frac{x_1-x_0}{2}+x_1-x_2),\no\\
&&\sigma_{3}=(1,2):\;\;\;\; x_0x_1x_2\frac{x_1+x_2}{2}=x_0x_2(\frac{2}{9}x_0^2+\frac{5}{9}x_0x_2+\frac{2}{9}x_2^2+
r_2(\frac{2}{9}x_0+\frac{1}{3}x_1+\frac{4}{9}x_2)),\no\\
&&(r_2=x_1-x_0+\frac{x_1-x_2}{2}),\no\\
&&\sigma_{3}=(1,1,1):\;\;\;\; x_0x_1x_2x_3=x_0x_3(\frac{2}{9}x_0^2+\frac{5}{9}x_0x_3+\frac{2}{9}x_3^2+
r_3(\frac{2}{9}x_0+\frac{1}{3}x_1+\frac{4}{9}x_3)+\no\\
&&+r_4(\frac{4}{9}x_0+\frac{1}{3}x_2+\frac{2}{9}x_3)+\frac{1}{3}r_3r_4)
,\no\\
&&(r_3=2x_1-x_0-x_2,\;r_4=2x_2-x_1-x_3).
\end{eqnarray}
With these decompositions, the same argument on residue integrals as the one used in the proof of Theorem 2 leads us 
to the desired recursive formulas. We can prove the recursive formulas for higher degree by extending this kind of discussion.

\end{document}